# CORESTRICTION PRINCIPLE IN NON-ABELIAN GALOIS COHOMOLOGY

Nguyêñ Quôć Thăńg *

Department of Mathematics, Israel Institute of Technology
Technion, Haifa - 32000
Israel

**Abstract**

In this paper we prove that over local or global fields of characteristic 0, the Corestriction Principle holds for kernel and image of all maps which are connecting maps in group cohomology and the groups of $R$-equivalences. Some related questions over arbitrary fields of characteristic 0 are also discussed.

**Plan.**

## Introduction.

*Supported by a McMaster University and Lady Davis Fellowships. E-mail : nguyen@tx.technion.ac.il

Let $k$ be a field of characteristic 0, $G$ a linear algebraic group defined over $k$. We are interested only in linear $k$-groups, so the adjective "linear" is omitted.

It is well-known that if $G$ is commutative, then for any finite extension $k'$ of $k$, there is the so-called corestriction map

$$Cores_G : \mathrm{H}^q(k', G) \to \mathrm{H}^q(k, G), q \geq 0,$$

where $\mathrm{H}^q(L, H)$ denotes the Galois cohomology $\mathrm{H}^q(Gal(\bar{L}/L), H(\bar{L}))$ for a $L$-group $H$ defined over a field $L$ of characteristic 0 (or a perfect field $L$).

However if $G$ is not commutative, there is no such a map in general (see example 6) below), and, as far as we know, the most general sufficient conditions are given in [Ri1], under which such a map can be constructed. The Corestriction Theory constructed there has many applications to theory of algebras, representation theory and related questions (see also [Ri2]). In this paper we are interested in the following natural question about the corestriction map:

Assume that there is a map, which is functorial in $k$ :

$$\alpha : \mathrm{H}^p(k, G) \to \mathrm{H}^q(k, T),$$

where $T$ is a commutative $k$-group, $G$ a non-commutative $k$-group. By restriction, for any finite extension $k'/k$ we have a functorial map

$$\alpha' : \mathrm{H}^p(k', G) \to \mathrm{H}^q(k', T).$$

**Question.** *When does* $Cores_T(\mathrm{Im}\ (\alpha')) \subset \mathrm{Im}\ (\alpha)$?

Of course, if there exists $Cores_G$ (e.g. under the conditions given in [Ri1]), which is functorial then the above question always has an affirmative answer. If the answer is affirmative for all $k'$, we say that *the Corestriction Principle holds* for (the image of) the map $\alpha$. One defines similar notion for the kernel of a map $\beta : \mathrm{H}^p(k, T) \to \mathrm{H}^q(k, G)$.

We say that the map $\alpha : \mathrm{H}^p(k, G) \to \mathrm{H}^q(k, T)$ is *standard* if it is obtained as a connecting map from the exact cohomology sequence associated with an exact sequence of $k$-groups involving $G$ and $T$. For example, let

$$1 \to A \to B \to C \to 1,$$

be an exact sequence of $k$-groups, where $A$ is considered as a normal $k$-subgroup of $B$. Then

$$\mathrm{H}^i(k, A) \to \mathrm{H}^i(k, B), i = 0, 1,$$

and

$$\mathrm{H}^0(k, C) \to \mathrm{H}^1(k, A)$$

are standard maps. In general, $C$ is just a quotient space and may not be a group. If $A$ is a central subgroup of $G$, then $C$ is a group, and one may define a connecting standard map $\mathrm{H}^1(k, C) \to \mathrm{H}^2(k, A)$.

It is worth mentioning that in some particular cases, the above question has an affirmative answer unconditionally and the *Norm Principle* is said to hold if it holds for $p = q = 0$ (which approves the adjective *norm*).

*Examples.* 1) Let $D$ be a finite dimensional central simple algebra over $k$, $G$ the $k$-group defined by the condition $G(k) = \mathrm{GL}_n(D)$ (a $k$-form of the general linear group), $G' = [G, G]$ (the group defined by the condition $G'(k) = \mathrm{SL}_n(D)$). We have the following exact sequence of $k$-groups

$$1 \to G' \to G \stackrel{N}{\to} \mathbf{G}_m \to 1,$$

where $N$ denotes the map induced from the reduced norm $\mathrm{GL}_n(D) \stackrel{Nrd}{\to} k^*$.

It is well-known that

$$N_{k'/k}(Nrd((D \otimes k')^*)) \subset Nrd_{D/k}(D^*),$$

which says that the Corestriction Principle holds for the image of $\alpha = N$, $p = q = 0$.

2) Let $\Phi$ be a non-degenerate $J$-hermitian form with values in a division $k$-algebra $D$ of center $k_0$, which is $k$ (resp. a separable quadratic extension of $k$), if the involution $J$ of $D$ is of the first (resp. second kind). Let $\mathrm{U}(\Phi)$ resp. $\mathrm{GU}(\Phi)$ be the $k$-group defined by the unitary group (resp. by the group of similarities) of the form $\Phi$. We have the following exact sequence of $k$-groups

$$1 \to \mathrm{U}(\Phi) \to \mathrm{GU}(\Phi) \stackrel{m}{\to} \mathbf{G}_m \to 1,$$

where the map $m$ maps every similarity to its similarity factor. It is known (see [L], [Sc] for the case of quadratic forms and [T1] for the case of skew-hermitian forms) that the Scharlau Norm Principle holds for the group of similarity factors, so the Corestriction Principle holds for the image of $\alpha = m$ and $p = q = 0$. Notice also that since $\mathrm{SU}(\Phi)$ is the connected component of $\mathrm{U}(\Phi)$ in the Zariski topology, it follows that the Norm Principle also holds for the group of special (or proper) similarity factors. We have

**Theorem.** *Let $\Phi$ be a non-degenerate skew-hermitian with values in a division $k$-algebra $D$ with respect to an involution $J$ of the first kind of $D$, trivial on $k$, and $M(\Phi)$ (resp. $M(\Phi)^+$) be the group of similarity factors of similitudes (resp. proper similitudes ) of $\Phi$. Then for any finite extension $k'$ of $k$,*

$$N_{k'/k}(M(\Phi \otimes k')) \subset M(\Phi),$$
$$N_{k'/k}(M(\Phi \otimes k')^+) \subset M(\Phi)^+.$$

3) Let $f$ be a non-degenerate quadratic form over a field $k$ of characteristic $\neq 2$. Let $\mathrm{Spin}(f)$ (resp. $\mathrm{SO}(f)$ be the Spin (resp. special orthogonal) $k$-group of $f$. Let $\mu_2$ be the group $\{\pm 1\}$. We have the following exact sequence

$$1 \to \mu_2 \to \mathrm{Spin}(f) \to \mathrm{SO}(f) \to 1,$$

hence also

$$\mathrm{Spin}(f)(k) \to \mathrm{SO}(f)(k) \xrightarrow{\delta} k^*/k^{*2}.$$

The Knebusch Norm Principle (see [L], or Section 3 below) allows one to deduce the Corestriction Principle for the image of $\delta$, $p = 0, q = 1$, which means that the Norm Principle holds for the spinor norms.

4) A new kind of Corestriction Principle over local and global fields has been found by P. Deligne [De, Prop. 2.4.8], which, in the case of characteristic 0 and in notations of abelian Galois cohomology ([B1], [Mi, Appendix B]), says that the Corestriction Principle for images holds for the map

$$ab^0_G : \mathrm{H}^0(k, G) \to \mathrm{H}^0_{ab}(k, G).$$

This result has been subsequently applied to various problems related with canonical models of Shimura varieties.

5) There are few other examples due to Gille [G1] and Merkurjev [M1] (see also Section 3 below), who proved that the Corestriction Principle holds for the image when restricting $\alpha$ to the *subgroup* $RG(k)$ of elements of $G(k)$ which are R-equivalent to 1.

6) Given any natural numbers $n \geq 2$, $r \geq 1$, Rosset and Tate have constructed in [RT] an example of a field $E$ containing the group $\mu_n$ of $n$-th roots of 1, a finite Galois extension $F$ of $E$ of degree $r$, and an element $x$ of $K_2(F)$, which is a symbol, such that the image of $x$ via the trace

$$Tr_{F/E} : K_2F \to K_2E$$

is a sum of at least $r$ symbols. From this they derive a symbol algebra of degree $n$ over $F$, considered as an element of $\mathrm{H}^2(F, \mu_n)$, such that its image via the corestriction

$$Cores_{F/E} : \mathrm{H}^2(F, \mu_n) \to \mathrm{H}^2(E, \mu_n)$$

is not a symbol. Therefore the question above has a negative answer for the standard map

$$\Delta : \mathrm{H}^1(E, \mathrm{PGL}_n) \to \mathrm{H}^2(E, \mu_n).$$

Despite of this, we will see that in many interesting cases, the Corestriction Principle for standard maps hold. The purpose of this paper is to discuss the validity of the Corestriction Principle for images and kernels of standard maps in the case the field of definition is a local or global field of characteristic 0, and its applications. If the base field is an arbitrary field of characteristic 0, we discuss the relation between the corestriction principles for various types of standard maps. As applications, we give a new proof of Merkurjev's Norm Principle and prove the Corestriction Principle for (the images of) maps $\pi_R : G(k)/R \to T(k)/R$, where $G, T$ are connected reductive groups with $T$ commutative, $G(k)/R$ and $T(k)/R$ denote the corresponding groups of $R$-equivalences and $\pi_R$ is induced from a $k$-homomorphism $\pi : G \to T$. The reason that we insist on calling *corestriction principle* is that indeed, all the resulting "norm maps" are induced from certain corestriction maps in usual cohomology theory.

# 1 Corestriction Principle in non-abelian cohomology : local and global fields.

In this section we prove the validity of the Corestriction Principle for images and kernels of standard maps for local or global base fields of characteristic 0 and consider some applications.

Our first main result of this section is the following

**1.1. Theorem.** *Let $k$ be a local or global field of characteristic 0, $G$ a connected $k$-group, $T$ a connected commutative $k$-group and $\alpha : \mathrm{H}^p(k, G) \to \mathrm{H}^q(k, T)$ a standard map. Assume that $G$ is a central extension of $T$ if $p = q = 2$ where the 2-cohomology is defined as in [Gi]. Then for $0 \leq p \leq q \leq 2$ the Corestriction Principle holds for the image of $\alpha$.*

*Proof.* We assume the familiarity with the notion and results from the Borovoi - Kottwitz theory of abelian Galois cohomology of algebraic groups as presented in [B1-3] (see also [Mi, Appendix B] for a survey). We may assume also that $G$ is not abelian.

We first begin with the case of small $p, q$.

a) Let $p = q = 0$. Then we may assume that the map (denoted by the same symbol) $\alpha : G \to T$, induced from $\alpha : G(k) \to T(k)$, is surjective. Then we have the following exact sequence of $k$-groups

$$1 \to G_1 \to G \stackrel{\alpha}{\to} T \to 1,$$

with $G_1 = \mathrm{Ker}\ (\alpha)$. It is easy to see that $\alpha$ is surjective on $R_u(G)(k)$, i.e., $\alpha(R_u(G)(k)) = R_u(T)(k)$, where $R_u(.)$ denotes the unipotent radical of (.). Hence we may assume that $G$ is reductive and $T$ is a torus. Therefore $G_1$ contains $G' = [G, G]$.

Let $G = G'T'$, $F = \mathrm{Ker}\ (\tilde{G} \stackrel{\rho}{\to} G')$, where $\tilde{G}$ denotes the simply connected covering of $G'$. First we assume that $G_1 = G'$. By Proposition 2.4.8 of [De], there exists a corestriction map

$$Cores : G(k')/\rho(\tilde{G}(k')) \to G(k)/\rho(\tilde{G}(k)).$$

(The proof of Deligne [De] and [B1], [B3] show that in fact Deligne has proved the Corestriction Principle for $ab^0$ for any connected reductive group over local or global fields of characteristic 0.) We claim that this map, while restricted to a subgroup $H(k')/\rho(\tilde{G}(k'))$, where $H$ is a connected $k$-subgroup of $G$, containing $G'$, is the one constructed by Deligne.

Indeed, we have the following commutative diagram

$$\begin{array}{ccc} \mathrm{H}^0(k',H) & \to & \mathrm{H}^0_{ab}(k',H) \\ \downarrow & & \downarrow \\ \mathrm{H}^0(k',G) & \to & \mathrm{H}^0_{ab}(k',G), \end{array}$$

where all maps are functorial (see [B1]). Then the image of $\mathrm{H}^0(k',H)$ in $\mathrm{H}^0_{ab}(k',H)$ is $H(k')/\rho(\tilde{G}(k'))$ by [B1]. Therefore the claim follows when we project this diagram into similar diagram where $k'$ is replaced by $k$ and by making use of the commutativity of suitable related diagrams. (We can state in fact a more general statement, but we do not need it here.)

Thus we have the following *commutative* diagram with exact rows

$$\begin{array}{ccccccc} 1 \to & G'(k')/\rho(\tilde{G}(k')) & \to & G(k')/\rho(\tilde{G}(k')) & \to & \mathrm{H}^0(k',T) \\ & \downarrow & & \downarrow & & \downarrow \\ 1 \to & G'(k)/\rho(\tilde{G}(k)) & \to & G(k)/\rho(\tilde{G}(k)) & \to & \mathrm{H}^0(k,T), \end{array}$$

hence also the following corestriction (norm) map

$$G(k')/G'(k') \to G(k)/G'(k). \tag{1}$$

Since these two groups are respectively the images of $G(k')$ and $G(k)$ in $\mathrm{H}^0(k',T)$ and $\mathrm{H}^0(k,T)$, the assertion of the theorem follows.

Now we turn to the general case. Let us consider the following commutative diagram

$$\begin{array}{ccccccccc} 1 \to & G' & \to & G & \to & T' & \to 1 \\ & \downarrow & & \downarrow & & \downarrow & \\ 1 \to & G_1 & \to & G & \to & T & \to 1, \end{array}$$

where $G_1$ is any $k$-subgroup of $G$ containing $G' = [G,G]$. By taking the induced commutative diagram of exact cohomology sequences of these two rows and by using the fact that the Corestriction Principle holds for the image of the map $G(k) \to T'(k)$ shown above, we obtain the Corestriction Principle for the image of $G(k) \to T(k)$.

b) $p = 0, q = 1$. Let
$$1 \to T \to G_1 \to G \to 1$$
be the exact sequence of $k$-groups under consideration. Since $T = T_s \times T_u$, where $T_s$ (resp. $T_u$) is a $k$-torus (resp. unipotent $k$-group), $\mathrm{H}^1(k, T_u) = 0$ and $1 \to T_u \to R_u(G_1) \to R_u(G) \to 1$ is an exact sequence, we may assume that $G_1, G$ and $T$ are reductive.

Then by [B1] we have the following commutative diagram

$$\begin{array}{ccccc} \mathrm{H}^0(k, G_1) & \to & \mathrm{H}^0(k, G) & \xrightarrow{\alpha} & \mathrm{H}^1(k, T) \\ \downarrow & & \downarrow & & \downarrow \\ \mathrm{H}^0_{ab}(k, G_1) & \to & \mathrm{H}^0_{ab}(k, G) & \to & \mathrm{H}^1_{ab}(k, T), \end{array}$$

where $\mathrm{H}^i_{ab}(.,.)$ denotes the $i$-th abelian cohomology and the vertical maps are the maps $ab^i$ constructed in [B1]. Since $\mathrm{H}^1(k, T) \simeq \mathrm{H}^1_{ab}(k, T)$, it follows that if $ab^0_G$ satisfies the Corestriction Principle (for images), then $\alpha$ does also. By [B1, p.39, Proposition 3.6] we have

$$\mathrm{Im}\ (\mathrm{H}^0(k, G) \to \mathrm{H}^0_{ab}(k, G)) = G(k)/\rho(\tilde{G}(k)),$$

hence by making use of the Deligne map above the assertion is true in this case.

c) $p = q = 1$ or $p = 1, q = 2$. It is well-known that there is a canonical bijection between $\mathrm{H}^1(k, G)$ and $\mathrm{H}^1(k, L)$, where $L$ is any Levi $k$-subgroup of $G$ and that $\mathrm{H}^2(k, U) = 0$ for any commutative unipotent $k$-group by a theorem of Serre. Hence we may assume that $G$ and $T$ are reductive. We have the following commutative diagram

$$\begin{array}{ccc} \mathrm{H}^1(k,G) & \to & \mathrm{H}^1_{ab}(k,G) \\ \downarrow & & \downarrow \\ \mathrm{H}^q(k,T) & \to & \mathrm{H}^q_{ab}(k,T) \end{array}$$

where $q = 1, 2$ (see [B1]). Since $ab^1_G$ is surjective for local or global fields of charcteristic 0 ([B1]) and since $\mathrm{H}^1_{ab}(k', G) \to \mathrm{H}^1_{ab}(k, G)$ exists, the assertion of the theorem is verified.

d) $p = q = 2$. Let $1 \to G_1 \to G \to T \to 1$ be the exact sequence of $k$-groups under consideration. We have the following exact sequence of cohomology (by assumption)

$$\mathrm{H}^1(k,T) \to \mathrm{H}^2(k,G_1) \to \mathrm{H}^2(k,G) \to \mathrm{H}^2(k,T),$$

where $G_1$ and $T$ are commutative and $G_1$ is central subgroup of $G$. It follows that $G$ is solvable. If $T = T_s \times T_u$, where $T_u$ is the unipotent part of $T$, then we know that $\mathrm{H}^2(k, T_u) = 0$, hence we may assume that $T$ is a torus. Then $G_1$ contains $R_u(G)$ so $G$ is a connected nilpotent group, for which the assertion is obvious. ∎

**1.2. Remarks.** 1) It follows from the construction of $ab^2_G$ of [B2, p. 228] that this map satisfies the Corestriction Principle for images for any field $k$ of characteristic 0 and any connected reductive $k$-group $G$.

2) It is desirable to modify the Borovoi - Kottwitz theory so that it can cover also the case where the characteristic of $k$ is $p > 0$.

3) One may define the corestriction map (or "norm map") between some factor sets of $\mathrm{H}^i(k, G)$, e.g., in the following cases ($k$ is a local or global field):
$p = 0, q \leq 1$. Then we obtain indeed a norm map

$$N : Coker(\alpha \otimes k') \to Coker(\alpha).$$

$p = 1, q = 2$, $T$ is a central subgroup of $G$.

To be complete, together with the Corestriction Principle for the *images* of standard maps, we need also to consider the validity of this principle for *kernels* of standard maps. Namely for a standard map

$$\alpha : \mathrm{H}^p(k, T) \to \mathrm{H}^q(k, G),$$

where $T$, $G$ are connected $k$-groups with $T$ commutative, and for a finite extension $k'$ of $k$ with the corestriction map $Cores_T : \mathrm{H}^p(k', T) \to \mathrm{H}^p(k, T)$, we ask

**Question.** *When does* $Cores_T(\mathrm{Ker}\ (\alpha \otimes k')) \subset \mathrm{Ker}\ (\alpha)$ ?

By using Theorem 1.1 it is easy to see that in the case $k$ is a global or a global field of characteristic 0, one is reduced to considering the case $p = q = 1$. We have the following affirmative result for local and global fields of characteristic 0.

**1.3. Theorem.** *Let $k$ be a local or global field of characteristic 0 and $T$ a connected commutative $k$-subgroup of a connected $k$-group $G$. Then the Corestriction Principle holds for the kernel of the standard map $\alpha : \mathrm{H}^1(k, T) \to \mathrm{H}^1(k, G)$.*

*Proof.* As above, we may assume that $T$ is a $k$-subtorus of $G$ and $G$ is reductive. We need the following lemmas.

**1.4. Lemma.** *Assume that we have the following commutative diagram*

$$\begin{array}{ccccc} A' & \xrightarrow{p'} & B' & \xrightarrow{q'} & C' \\ \downarrow \beta & & \downarrow \gamma & & \\ A & \xrightarrow{p} & B & \xrightarrow{q} & C, \end{array}$$

*where $A$, $B$, $A'$, $B'$ are groups, the left diagram is a commutative diagram of groups. Let $e' = q'(1)$, $e = q(1)$, where 1 is the identity element of the corresponding groups. Then if $\gamma(q'^{-1}(e')) \subset q^{-1}(e)$ then*

$$\beta(r'^{-1}(e')) \subset r^{-1}(e),$$

*with* $r' = q'p'$, $r = qp$.

*Proof.* We have

$$\begin{aligned} r(\beta(r'^{-1}(e'))) &= q(p(\beta(r'^{-1}(e')))) \\ &= q(\gamma(p'(r'^{-1}(e')))) \\ &= q(\gamma(p'(p'^{-1}(q'^{-1}(e'))))) \\ &\subset q(\gamma(q'^{-1}(e'))) \\ &\subset q(q'^{-1}(e)) = e \end{aligned}$$

and the lemma follows. ∎

Recall that a connected reductive $k$-group $H$ is a $z$-extension of a $k$-group $G$ if $H$ is an extension of $G$ by an *induced* $k$-torus $Z$, such that the derived subgroup (semisimple part) $[H, H]$ of $H$ is simply connected. For a field extension $K/k$ and an element $x \in \mathrm{H}^1(K, G)$, a $z$-extension of $G$ over $k$ is called $x$-lifting if $x \in \mathrm{Im}\ (\mathrm{H}^1(K, H) \to \mathrm{H}^1(K, G))$.

**1.5. Lemma.** *Let $G$ be a connected reductive $k$-group, $K$ a finite extension of $k$, $x$ an element of $\mathrm{H}^1(K, G)$. Then there is a $z$-extension*

$$1 \to Z \to H \to G \to 1,$$

*of $G$, where all groups and morphisms are defined over $k$, which is $x$-lifting.*

**1.6. Lemma.** *Let $\alpha : G_1 \to G_2$ be a homomorphism of connected reductive groups, all defined over $k$, $x \in \mathrm{H}^1(K, G_1)$, where $K$ is a finite extension of $k$. Then there exists a $x$-lifting $z$-extension $\alpha' : H_1 \to H_2$ of $\alpha$, i.e., $H_i$ is a $z$-extension of $G_i$ $(i = 1, 2)$, and we have the following commutative diagram*

$$\begin{array}{ccc} H_1 & \xrightarrow{\alpha'} & H_2 \\ \downarrow & & \downarrow \\ G_1 & \xrightarrow{\alpha} & G_2, \end{array}$$

*with all groups and morphisms defined over $k$.*

The Lemmas 1.5 - 1.6, in the case $K = k$, are due to Kottwitz (see e.g. [B1, p. 34 and p. 37]). The proofs in our case are the same : Lemma 1.6 follows from Lemma 1.5. To prove Lemma 1.5 we choose a Galois extension $F/k$ large enough so that $F$ contains $K$ and $x$ is split over $F$ (i.e. $res_{K/F}(x) = 1$, where $res_{K/F} : \mathrm{H}^1(K, G) \to \mathrm{H}^1(F, G)$), and such that there is a $z$-extension

$$1 \to Z \to H \to G \to 1$$

with

$$Z \simeq (\mathrm{R}_{F/k}(\mathbf{G}_m))^n$$

for some $n$ (see [B1, pp. 33 - 34] for more details). Then one checks that the image of $x$ in $\mathrm{H}^2(K, Z)$ is trivial. Hence $x \in \mathrm{Im}\ (\mathrm{H}^1(K, H) \to \mathrm{H}^1(K, G))$.

By Lemma 1.4, we may assume that $T$ is a maximal torus of $G$ and by Lemma 1.6, we may assume that $G$ has simply connected semisimple part.. In the case of local fields we give two arguments to prove the assertion of the theorem.

First, let $x \in \mathrm{Ker}\ (\alpha)$. By Lemma 1.6 there exists a $x$-lifting $z$-extension $T_1 \to G_1$ of $\alpha$, all defined over $k$. Since $T$ is a torus , $T_1$ is also a torus. It is easy to see that if the Corestriction Principle for kernels holds for any pair $(T_1, G_1)$ with $G_1$ having the simply connected semisimple part then it also holds for $(T, G)$. So from now on we assume that $G' = [G, G]$ is simply connected.

First we assume that $k$ is a local field. The case $k = \mathbf{R}$ is trivial, so we assume that $k$ is a $p$-adic field. Let $S$ be the maximal central torus of $G$, $G = SG'$. We have the following commutative diagram

$$
\begin{array}{ccc}
 & & \mathrm{H}^1(k,G') \\
 & & \downarrow p \\
\mathrm{H}^1(k,T) & \stackrel{\alpha}{\to} & \mathrm{H}^1(k,G) \\
 & & \downarrow q \\
 & & \mathrm{H}^1(k,G/G').
\end{array}
$$

Since $\mathrm{H}^1(k,G')=0$ by Kneser's Theorem, Ker $(q)=0$. Therefore

$$\mathrm{Ker}\ (\alpha)=\mathrm{Ker}\ (q\alpha).$$

Since $G/G'$ is a torus, $q\alpha : T \to G/G'$ satisfies the Corestriction Principle for kernels. Hence the assertion of the theorem is verified for local fields.

Now we assume that $k$ is either a local or a number field. By making use of the generalized Ono's trick due to Sansuc (see [Sa, Lemme 1.10]), we can find a natural number $m$, quasi-split (induced) $k$-tori $P,Q$ such that there is a central $k$-isogeny

$$1 \to F \to G_1 \to G^m \times Q \to 1,$$

where $F$ is a finite central subgroup of a connected reductive $k$-group $G_1$, which is a direct product of $P$ and a simply connected semisimple group $G'_1$. Let $T_1$ be the unique maximal $k$-torus of $G_1$ covering the maximal torus $T'=T^m\times Q$ of $G'=G^m\times Q$, $T_1=\tilde{T}\times P$, $G_1=\tilde{G}_1\times P$, where $\tilde{T}$ is a maximal torus of the semisimple simply connected (derived) subgroup $\tilde{G}_1=G'_1$ of $G_1$. It is clear that the assertion of the theorem for $(T,G)$ is equivalent to that for $(T',G')$. Recall that we may assume the semisimple part of $G'$ to be simply connected, i.e., isomorphic to $\tilde{G}_1$. Then $G'=\tilde{G}_1P'$, where $P'$ is the image of $P$. We have the following commutative diagram with exact rows

$$
\begin{array}{ccccccc}
\mathrm{H}^1(k,F) & \stackrel{\theta}{\to} & \mathrm{H}^1(k,T_1) & \stackrel{\beta}{\to} & \mathrm{H}^1(k,T') & \stackrel{\delta}{\to} & \mathrm{H}^2(k,F) \\
\downarrow = & & \downarrow \gamma & & \downarrow \alpha & & \downarrow = \\
\mathrm{H}^1(k,F) & \stackrel{p}{\to} & \mathrm{H}^1(k,G_1) & \stackrel{\pi}{\to} & \mathrm{H}^1(k,G') & \stackrel{\delta}{\to} & \mathrm{H}^2(k,F)
\end{array}
$$

Note that if $\alpha(x) = 0$, then $0 = \delta(\alpha(x)) = \delta(x)$, hence $x \in \mathrm{Im}\ (\beta)$ and

$$\mathrm{Ker}\ (\alpha) = \beta\ (\mathrm{Ker}\ (\alpha\beta))$$

$$= \beta(\{x \in \mathrm{H}^1(k, T_1) : \gamma(x) \in \mathrm{Im}\ (p)\}).$$

Hence it suffices to show that for the set

$$A(k) := \{x \in \mathrm{H}^1(k, T_1) : \gamma(x) \in \mathrm{Im}\ (p)\},$$

we have

$$Cores_{k'/k}(A(k')) \subset A(k).$$

Since $P$ is an induced torus, we have $\mathrm{H}^1(k, P) = 0$, and

$$\mathrm{H}^1(k, T_1) = \mathrm{H}^1(k, \tilde{T_1}) \times \{0\}, \mathrm{H}^1(k, G_1) = \mathrm{H}^1(k, \tilde{G_1}) \times \{0\},$$

hence $A(k)$ may be identified with the following set $\{x \in \mathrm{H}^1(k, \tilde{T_1}) : \pi(\gamma(x)) = 0\}$, where $\pi$ may be considered as the map, induced from the composition

$$\tilde{G}_1 \hookrightarrow G' = \tilde{G}_1 P',$$

since $\tilde{G}_1$ is simply connected so the restriction of $\pi$ on $\tilde{G}_1$ is an isomorphism. Let $F' := \tilde{G}_1 \cap P'$, $\bar{P} = P'/F'$. Then we have the following commutative diagram with exact rows

$$\begin{array}{ccccc} & & \mathrm{H}^1(k, \tilde{T_1}) & \xrightarrow{=} & \mathrm{H}^1(k, \tilde{T_1}) \\ & & \downarrow i & & \downarrow \\ \bar{S}(k) & \xrightarrow{\delta_k} & \mathrm{H}^1(k, \tilde{G}_1) & \rightarrow & \mathrm{H}^1(k, \tilde{G}_1 P') \end{array}$$

Let $t' \in A(k')$ and $(t'_s)$ be a representative of $t'$, $t'_s \in T_1(k_s)$. Then $t'_s = (g'p')^{-1}\ {}^s(g'p')$ for some $g' \in \tilde{G}_1(k_s)$ and $p' \in P'(k_s)$ and for all $s \in Gal(k_s/k')$. It follows that $f'_s := p'^{-1}\ {}^s p' \in F'(k_s)$ is a cocycle, representing an element $f'$ from Ker $(\mathrm{H}^1(k', F') \to \mathrm{H}^1(k', P')$ and we see that $t''_s := t'_s f'^{-1}_s = g'^{-1}\ {}^s g'$ represents an element $t''$ from Ker $(\mathrm{H}^1(k', \tilde{T_1}) \to \mathrm{H}^1(k', \tilde{G}_1)$. Then

$$Cores_{k'/k}(t'') = Cores_{k'/k}(t') Cores_{k'/k}(f'^{-1}).$$

Since

$$Cores_{k'/k}(\text{Ker }(\mathrm{H}^1(k',F') \to \mathrm{H}^1(k',P')) \subset \text{Ker }(\mathrm{H}^1(k,F') \to \mathrm{H}^1(k,P')),$$

we can choose a representative $(f_r)_r$ of $Cores_{k'/k}(f')$, $r \in Gal(k_s/k)$, such that

$$f_r = p^{-1}\ {}^r p, \forall r \in Gal(k_s/k).$$

Assume that

$$Cores_{k'/k}(t'') = t \in \text{Ker }(\mathrm{H}^1(k,\tilde{T}_1) \to \mathrm{H}^1(k,\tilde{G}_1)),$$

$t = [(t_r)],\ t_r = g^{-1}\ {}^r g$, $g \in \tilde{G}_1(k_s)$, $r \in Gal(k_s/k)$. Then

$$Cores_{k'/k}(t') = Cores_{k'/k}(t'')Cores_{k'/k}(f')$$

has a representative $(t_{1,r})_r$, where

$$t_{1,r} = g^{-1}\ {}^r g f_r = g^{-1}\ {}^r g p^{-1}\ {}^r p = (gp)^{-1}\ {}^r(gp),$$

i.e., $Cores_{k'/k}(t') \in \text{Ker }(\mathrm{H}^1(k,\tilde{T}_1) \to \mathrm{H}^1(k,G'))$ as required.

Therefore we are reduced to proving the Corestriction Principle for kernels for $(\tilde{T}_1,\tilde{G}_1)$.

If $k$ is a $p$-adic field, then the assertion now is trivial due to the fact that $\mathrm{H}^1(k,\tilde{G}_1) = 0$. If $k$ is a number field, we have the following commutative diagram

$$\mathrm{H}^1(k,\tilde{T}_1) \stackrel{\alpha_1}{\to} \mathrm{H}^1(k,\tilde{G}_1)$$

$$\downarrow \lambda \qquad\qquad \downarrow \lambda'$$

$$\textstyle\prod_{v\in\infty} \mathrm{H}^1(k_v,\tilde{T}_1) \stackrel{\alpha'_1}{\to} \prod_{v\in\infty} \mathrm{H}^1(k_v,\tilde{G}_1),$$

where $\infty$ denotes the set of infinite places of $k$. Since the cohomological Hasse principle holds for $\mathrm{H}^1$ of simply connected semisimple $k$-groups, Ker $(\lambda') = 0$, hence

$$\text{Ker }(\alpha_1) = \text{Ker }(\lambda'\alpha_1) = \text{Ker }(\alpha'_1\lambda).$$

By Lemma 1.4 and the local field case above the proof of Theorem 1.3 follows from the last equality. ■

From the proof of Theorem 1.3 and results of Section 2 we derive the following.

**1.7. Corollary.** *The Corestriction Principle for kernels of the standard maps* $\mathrm{H}^1(k,T) \to \mathrm{H}^1(k,G)$*, where* $T$ *and* $G$ *are connected groups over a field* $k$ *of characteristic 0,* $T$ *is commutative, holds if and only if the same holds for all pairs* $(T,G)$ *with* $T$ *a maximal torus of a simply connected almost simple* $k$*-group* $G$*, all defined over* $k$*.*

**1.8. Remarks.** 1) The Corestriction Principle for kernels suggests the study of the kernels of maps $\mathrm{H}^1(k,T) \to \mathrm{H}^1(k,G)$, which are little known except for the case of local or global fields. It worth noticing that the study of such kernels plays an important role in the proof of the Hasse principle for $\mathrm{H}^1$ of simply connected semsimple groups, done by Harder ([Ha]). See also further comments done by Tits [Ti]. Moreover, the proof above shows that if $\mathrm{H}^1(k,\tilde{G}) = 0$, where $\tilde{G}$ is the semisimple simply connected covering of $G'$ (e.g., according to Bruhat - Tits, when $k$ is a local field with residue field of cohomological dimension $\leq 1$), the Corestriction Principle for kernels for $\mathrm{H}^1(k,T) \to \mathrm{H}^1(k,G)$ holds.

2) It is easy to show that for connected reductive groups $G$ over number fields $k$ there are norm maps $A(k',G) \to A(k,G)$ and $\mathrm{III}(k',G) \to \mathrm{III}(k,G)$ for all finite extensions $k \subset k'$, where $A(K,G)$ denotes the (defect of weak approximation) quotient group $\prod_v G(K_v)/Cl(G(K))$, where $Cl$ denotes the closure in the product topology of $G(K_v)$, and $\mathrm{III}(K,G)$ denotes the Tate - Shafarevich group of $G$. The first follows from a result of Sansuc [Sa, Thm 3.3], and the second follows from a result of Borovoi [B1, Thm 5.13].

3) It might be of interest to investigate the Norm Principle for the map $\alpha : X(k) \to T(k)$, where $X(k)$, $T(k)$ are some "objects over $k$" and $T(k)$ is a commutative group. More precisely, the example we have in mind is the following.

Given a non-constant $k$-rational map $\phi : X \to T$ from an irreducible $k$-variety $X$ into a commutative $k$-group $T$. One asks when $N(\phi(X(k')) \subset \phi(X(k))$. This and related questions will be the subject of a future study.

4) Let $G$ be a connected reductive group over a field $k$. As in the case of

semisimple groups, we define the *Whitehead group of $G$* over $k$, $W(k, G) := G(k)/G(k)^+$, where $G(k)^+$ denotes the subgroup of $G(k)$ generated by $k$-rational points of unipotent radicals of parabolic $k$-subgroups of $G$. Note that $G(k)^+$ is a normal subgroup of $G(k)$. It is known that over any local field (resp. global field) $k$, the Kneser - Tits conjecture holds for all isotropic simply connected almost simple groups $H$ (resp. except possibly for some groups of type ${}^2\mathrm{E}_6$) over $k$, i.e., $H(k) = H(k)^+$. Thus the Deligne 's norm map gives rise to the norm map for the Whitehead groups of connected reductive groups with isotropic almost simple factors (containing no almost simple factors of type ${}^2\mathrm{E}_6$ if $k$ is a number field). In particular the following natural question arises :

**Question.** *Let $k$ be an infinite field and $G$ be a connected reductive $k$-group. Is there any "norm relation" between $W(k', G)$ and $W(k, G)$ for all finite extension $k \subset k'$ ?*

For the case of a local or global field we will give an answer to this question in a relative form, namely modulo the image of the Whitehead group of a connected reductive $k$-group $G_0$ with semisimple part isogeneous to that of $G$ (see the corollary below).

5) In the case $p = q = 0$ we have seen that there are corestriction (norm) maps for the following quotient groups of $G(k)$ : $G(k)/G'(k)$ and $G(k)/\rho(\tilde{G}(k))$. It is natural to ask if there is similar map for other "intermediate" quotient groups, namely for $G(k)/\pi(G_0(k))$, where $G_0$ is a connected reductive $k$-group with a $k$-homomorphism $\pi : G_0 \to G$, which restricted to $G'_0$ is an isogeny onto the semisimple part $G'$ of $G$. The answer is affirmative and we have the following result, which is a slight generalization of a result of Deligne [De, Proposition 2.4.8].

**1.9. Theorem.** *With the above notation, assume that $G$ is a connected reductive $k$-group. For any finite extension $k'$ of a local or global field $k$ there is a canonical norm map*

$$G(k')/\pi(G_0(k')) \to G(k)/\pi(G_0(k)).$$

By taking the resctricted product of all such maps in the local case, as in

[De, 2.4.9] we deduce from the theorem the following

**1.10. Corollary.** *With above notation, we have a norm map*

$$N_{k'/k} : G(\mathbf{A}')/\pi(G_0(\mathbf{A}')) \to G(\mathbf{A})/\pi(G_0(\mathbf{A}),$$

*where* $\mathbf{A}'$, $\mathbf{A}$ *denotes the adele ring of* $k'$, $k$ *, respectively.*

*Proof of Theorem 1.9.* First we need the following

**1.11. Lemma.** *There is a canonical norm map*

$$G'(k')/\pi(G'_0(k')) \to G'(k)/\pi(G'_0(k)).$$

To see this, we consider the following commutative diagram

$$\begin{array}{ccccc} \tilde{G}(k') & \xrightarrow{p} & G'(k') & \xrightarrow{\delta'} & \mathrm{H}^1(k', F) \\ \downarrow & & \downarrow & & \downarrow \gamma \\ G'_0(k') & \xrightarrow{\pi} & G'(k') & \xrightarrow{\delta} & \mathrm{H}^1(k', B), \end{array}$$

where $F = \mathrm{Ker}\,(\tilde{G} \to G)$ and $B = \mathrm{Ker}\,(\tilde{G} \to G'_0)$. (Recall that $\tilde{G}$ is the simply connected covering for both $G'_0$ and $G'$.) Now $G'(k')/\pi(G'_0(k'))$ is the image of $G'(k')$ in $\mathrm{H}^1(k', B)$ which is equal to $\gamma(\delta'(G'(k'))$. Since $\delta'$ and $\gamma$ satisfy the Corestriction Principle for images (see Theorem 1.1 and Corollary 2.10, Section 2 (below)), the same holds for $\delta$.

Now we come to the proof of the theorem. First we prove the theorem when $G_0 = G'_0$, i.e., the central torus part of $G_0$ is trivial. We have the following commutative diagram

$$\begin{array}{ccccccccc} 1 & \to & G'(k')/\rho(\tilde{G}(k')) & \to & G(k')/\rho(\tilde{G}(k')) & \to & G(k')/G'(k') & \to & 1 \\ & & \downarrow \alpha & & \downarrow \beta & & =\downarrow & & \\ 1 & \to & G'(k')/\pi(G'_0(k')) & \to & G(k')/\pi(G'_0(k')) & \to & G(k')/G'(k') & \to & 1. \end{array}$$

We denote by $\alpha', \beta', \gamma$ the corresponding canonical corestriction maps for $G'(k')/\rho(\tilde{G}(k'))$, $G(k')/\rho(\tilde{G}(k'))$ and $G(k')/G'(k')$, which exist by what we have proved above. Consider the following maps :

$$\alpha'' : G'(k)/\rho(\tilde{G}(k)) \to G'(k)/\pi(G_0(k)),$$

$$\beta'' : G(k)/\rho(\tilde{G}(k)) \to G(k)/\pi(G_0(k)).$$

Let $b \in G(k')/\pi(G_0(k'))$, $b_1 \in G(k')/\rho(\tilde{G}(k'))$ such that $\beta(b_1) = b$. It is natural to define the image of $b$ in $G(k)/\pi(G_0(k))$ by $\beta''(\beta'(b_1))$. If $\beta(b_2) = b$, then $b_1 = b_2 a$, $a \in \text{Ker }(\alpha) = \text{Ker }(\beta)$. Hence

$$\begin{aligned} \beta''(\beta'(b_1)) &= \beta''(\beta'(b_2 a)) \\ &= \beta''(\beta'(b_2)\beta'(a)). \end{aligned}$$

Since $\text{Ker }(\alpha) = \text{Ker }(\beta)$, one sees that $\beta'(\text{Ker }(\beta)) = \alpha'(\text{Ker }(\alpha)) \subset \text{Ker }(\alpha'') = \text{Ker }(\beta'')$. Thus

$$\beta''(\beta'(b_1)) = \beta''(\beta'(b_2))$$

as required.

In the general case, let $G_0 = G_0' S$, where $S$ is a central connected (torus) part of $G_0$. We have the following "conjectural" commutative diagram

$$\begin{array}{ccccccccc} 1 \to & \pi(G_0'(k'))/\pi(G_0'(k')) & \to & G(k')/\pi(G_0'(k')) & \to & G(k')/\pi(G_0(k')) & \to 1 \\ & \downarrow (?)\eta & & \downarrow \mu & & (?)\zeta \downarrow & \\ 1 \to & \pi(G_0(k))/\pi(G_0'(k)) & \to & G(k)/\pi(G_0'(k)) & \to & G(k)/\pi(G_0(k)) & \to 1. \end{array}$$

where (?) means a map to be proved existing. It is clear that $\zeta$ will exist if we can prove that $\eta$ exists. Thus we are reduced to proving the existence of the following conjectural commutative diagram

$$\begin{array}{ccccccccc} 1 \to & F(k')G_0'(k')/G_0'(k') & \to & G_0(k')/G_0'(k')) & \to & \pi(G_0(k'))/\pi(G_0'(k')) & \to 1 \\ & \downarrow (?)\theta & & \downarrow \epsilon & & (?)\kappa \downarrow & \\ 1 \to & F(k)G_0'(k)/G_0'(k) & \to & G_0(k)/G_0'(k) & \to & \pi(G_0(k))/\pi(G_0'(k)) & \to 1, \end{array}$$

thus also to the existence of $\theta$, since the existence of $\epsilon$ is known due to the proof of case $a$) of Theorem 1.1 (see (1)). Since for any extension $k \subset K$ we have $F(K)G'_0(K)/G'_0(K) = F(K)/F_0(K)$ where $F_0 := F \cap G'_0$ is a finite central $k$-subgroup of $G'_0$, $\theta$ is nothing else than the norm map induced from that of $F$ and $F_0$. ∎

From this theorem we deduce immediately the following

**1.12. Corollary.** *Let the notation be as above. Then $\pi$ induces a canonical norm homomorphism*

$$W(k', G)/\pi_*(W(k', G_0)) \to W(k, G)/\pi_*(W(k, G_0)),$$

*where $\pi_*$ denotes the homomorphism $W(., G_0) \to W(., G)$ induced from $\pi$.*

By combining Theorem 1.1 and Theorem 1.3 we have the following main result of this section.

**1.13. Theorem.** (Corestriction Principle) *Let $G$ be a connected and $T$ a commutative algebraic groups, all defined over a local or global field of characteristic* 0. *Assume that $\alpha_k : \mathrm{H}^p(k, G) \to \mathrm{H}^q(k, T)$ (resp. $\beta_k : \mathrm{H}^q(k, T) \to \mathrm{H}^p(k, G)$) be a standard map. Then for any finite extension $k'/k$ we have*

$$Cores_{k'/k}(\mathrm{Im}\ (\alpha_{k'})) \subset \mathrm{Im}\ (\alpha_k)$$

$$(resp.\ Cores_{k'/k}(\mathrm{Ker}\ (\beta_{k'})) \subset \mathrm{Ker}\ (\beta_k),$$

*where $Cores_{k'/k}$ denotes the corestriction map of cohomology of $T$.*

# 2 Corestriction Principle in non-abelian cohomology : arbitrary field of characteristic 0.

In this section we will discuss some relation between the validity of Corestriction Principles for standard maps of various type. As applications we apply the results obtained to give new proof of a result of Deligne that we used in Section 1 and of a result of Merkurjev about the Norm Principle for images of the set $RG(k)$ of elements $R$-equivalent to 1 of $G(k)$ (cf. Section 3). For simplicity we consider only reductive groups.

Let $k$ be a field of characteristic 0 and $\alpha : \mathrm{H}^p(k, G) \to \mathrm{H}^q(k, T)$ be a standard map, where $p = 0, 1$, $q \leq p + 1$, $G$ and $T$ are connected reductive $k$-groups, $T$ is a torus. Denote by $\tilde{G}$ (resp. $\bar{G}$) the simply connected covering (resp. the adjoint) group of the semisimple part of $G$, $\tilde{F} = \mathrm{Ker}\ (\tilde{G} \to \bar{G})$, $F' = \mathrm{Ker}\ (G' \to \bar{G})$, where $G'$ is the semisimple part of $G$. We consider the following statements.

*a) The Corestriction Principle for images holds for any such $\alpha$.*

*b) The Corestriction Principle for images holds for* $\mathrm{H}^p(k, \bar{G}) \to \mathrm{H}^{p+1}(k, F')$ *for $p = 0, 1$.*

*c) The Corestriction Principle for images holds for* $\mathrm{H}^p(k, \bar{G}) \to \mathrm{H}^{p+1}(k, \tilde{F})$, *for $p = 0, 1$.*

*d) The Corestriction Principle for images holds for $ab^p_G : \mathrm{H}^p(k, G) \to \mathrm{H}^p_{ab}(k, G)$ for any such $G$.*

We will show later that if one of these conditions holds (e.g. if $k$ is a local or global field) then for any isogeny of connected reductive k-groups $1 \to F \to G_1 \to G_2 \to 1$, the Corestriction Principle for the image of $\mathrm{H}^p(k, G_2) \to \mathrm{H}^{p+1}(k, F)$, $p = 0, 1$ holds.

We have the following results.

**2.1. Proposition.** *If d) holds then a) holds.*

*Proof.* The proof follows immediately from the functoriality of the maps $ab^p_G : \mathrm{H}^p(k, G) \to \mathrm{H}^p_{ab}(k, G)$, $p = 0, 1$, proved in [B1]. ∎

**2.2. Proposition.** *If d) holds for connected reductive k-groups with simply connected semisimple parts then d) holds ifself.*

*Proof.* For any finite extension $k'$ of $k$ let $\theta \in \mathrm{H}^p(k', G)$ be any element. We choose a $\theta$-lifting $z$-extension, all defined over $k :\ 1 \to Z \to H \to G \to 1$, which is possible due to Lemma 1.5. Recall that $H$ is a connected reductive $k$-group with simply connected semisimple part and $Z$ is a quasi-split $k$-torus. Let denote the induced (standard) maps

$$\phi : \mathrm{H}^p(k, H) \to \mathrm{H}^p(k, G),$$

$$\psi : \mathrm{H}^p(k', H) \to \mathrm{H}^p(k', G),$$

and let $\phi'$ and $\psi'$ stand for similar maps where $\mathrm{H}^p$ is replaced by $\mathrm{H}^p_{ab}$.

We have the following commutative diagram where all the vertical maps are the maps $ab^p$, where $ab'$ will denote the same map when we restrict to $k'$ :

$$\begin{array}{ccccccc}
 & & & & \mathrm{H}^p(k,H) & \xrightarrow{\psi} & \mathrm{H}^p(k,G) \\
\mathrm{H}^p(k',H) & \xrightarrow{\phi} & \mathrm{H}^p(k',G) & & \downarrow & & \downarrow \\
\downarrow & & \downarrow & & \mathrm{H}^p_{ab}(k,H) & \xrightarrow{\psi'} & \mathrm{H}^p_{ab}(k,G) \\
 & \nearrow & & \nearrow & & & \\
\mathrm{H}^p_{ab}(k',H) & \xrightarrow{\phi'} & \mathrm{H}^p_{ab}(k',G) & & & &
\end{array}$$

Let $\eta \in \mathrm{H}^p(k', H)$ such that $\phi(\eta) = \theta$. Then

$$ab'_G(\theta) = ab'_G(\phi(\eta)) = \phi'(ab_H(\eta)).$$

Assuming that $d)$ holds for $H$, there is $\alpha \in \mathrm{H}^p(k, H)$ such that

$$ab_H(\alpha) = Cor(ab'_H(\eta)).$$

Hence

$$\begin{aligned}
Cor(ab'_G(\theta)) &= Cor(ab'_G(\phi(\eta))) \\
&= Cor(\phi'(ab'_H(\eta))) \\
&= \psi'(Cor(ab'_H(\eta)))
\end{aligned}$$

$$= \psi'(ab_H(\alpha))$$

$$= ab_G(\psi(\alpha))$$

as required. ∎

We will need the following lemma in the sequel.

**2.3. Lemma.** *Assume that we are given the following diagram of pointed sets with distinguished elements*

$$\begin{array}{ccc} & & E \\ & & \downarrow \delta \\ & C \xrightarrow{\gamma} & D \\ & \downarrow \beta & \\ A \xrightarrow{\alpha} B & & \end{array}$$

*Then there exists a (non-uniquely determined) pointed set $F$ with morphisms (of pointed sets) $B \xrightarrow{p} F$ and $D \xrightarrow{q} F$ such that $p\beta = q\gamma$ and the following sequences are exact :*

$$A \xrightarrow{\alpha} B \xrightarrow{p} F,$$

$$E \xrightarrow{\delta} D \xrightarrow{q} F.$$

The proof of the lemma is trivial, so we omit it.

**2.4. Proposition.** *Assume that the Corestriction Principle for images holds for the map $\mathrm{H}^0(k,G) \to \mathrm{H}^0(k,T)$ (resp. for the map $\mathrm{H}^1(k,G) \to \mathrm{H}^1(k,T)$) for any $G$ and $T$ as above. Then the same holds for $ab_G^0 : \mathrm{H}^0(k,G) \to \mathrm{H}^0_{ab}(k,G)$ (resp. for $ab_G^1 : \mathrm{H}^1(k,G) \to \mathrm{H}^1_{ab}(k,G)$. In particular, if a) holds, then d) holds.*

*Proof.* By Proposition 2.2 we may assume that $G'$ is simply connected. By [B1] we have

$$\mathrm{H}^p_{ab}(k,G) = \mathrm{H}^p(k,G/G'),$$

hence $ab_G^p$ becomes just standard map ($p = 0, 1$). Since $G/G'$ is a torus, the proposition follows. ∎

**2.5. Proposition.** *We have b) $\Leftrightarrow$ c).*

*Proof.* We need only to prove that $c) \Rightarrow b)$. Consider the following commutative diagram.

$$\begin{array}{ccccc} \mathrm{H}^p(k', \tilde{G}) & \xrightarrow{p} & \mathrm{H}^p(k', \bar{G}) & \xrightarrow{\delta'} & \mathrm{H}^{p+1}(k', \tilde{F}) \\ \downarrow & & \downarrow & & \downarrow \gamma \\ \mathrm{H}^p(k', G') & \xrightarrow{\pi} & \mathrm{H}^p(k', \bar{G}) & \xrightarrow{\delta} & \mathrm{H}^{p+1}(k', F'), \end{array}$$

where $p = 0, 1$, $F' = \mathrm{Ker}\ (G \to \bar{G})$. (Recall that $\tilde{G}$ is the simply connected covering for both $\bar{G}$ and $G'$.) One sees that $\delta = \gamma\delta'$. Thus if the Corestriction Principle for images holds for $\delta'$, the same holds for $\delta$. ∎

**2.6. Proposition.** *Assume that a) holds for all $G$ with simply connected semisimple part $G'$. Then a) holds itself.*

*Proof.* For $p = 0$ it follows easily by taking any $z$-extension of $G$. For $p = q = 1$, it follows from Lemma 1.6 that for any finite extension $k'$ of $k$ and any element $x \in \mathrm{H}^1(k', G)$, there exists a $x$-lifting $z$-extension of $\pi : G \to T$, all defined over $k$ :

$$\begin{array}{ccc} H_1 & \xrightarrow{\pi'} & H_2 \\ \downarrow & & \downarrow \\ G & \xrightarrow{\pi} & T \end{array}$$

Here $H_2$ is a torus and $H_1$ has simply connected semisimple part. By assumption the Corestriction Principle holds for the image of $\mathrm{H}^1(k, H_1) \to \mathrm{H}^1(k, H_2)$. By chasing on suitable diagrams one sees that the image of $x$ in $\mathrm{H}^1(k, T)$ via $Cores_T\ \alpha$ lies in the image of $\mathrm{H}^1(k, G)$, where $\alpha : \mathrm{H}^1(k, G) \to \mathrm{H}^1(k, T)$. Hence the Corestriction Principle for images holds for $\alpha$.

The case $p = 1, q = 2$ is considered in a similar way. ∎

**2.7. Proposition.** *Assume that the Corestriction Principle for the image of* $\mathrm{H}^p(k, \bar{G}) \to \mathrm{H}^{p+1}(k, \tilde{F})$ *holds for all* $G$, $\tilde{F}$ *above, where p=0 (resp. p=1). Then the same holds for* $\mathrm{H}^p(k, G) \to \mathrm{H}^p(k, T)$ *for all* $G$, $T$ *above and for p=0 (resp. p=1). In particular, if c) holds then a) holds.*

*Proof.* By Proposition 2.6 we may assume that $G'$ is simply connected, $G' = \tilde{G}$.

Let $G = \tilde{G}.S$, where $S$ is a central torus of $G$, $F = \tilde{G} \cap S$ is a finite subgroup of $G$. First we consider the case $p = 0$.

Consider the following commutative diagram

$$\begin{array}{ccccccccc} 1 & \to & \tilde{F} & \to & \tilde{G}.S & \to & \bar{G} \times S/F & \to & 1 \\ & & \downarrow & & =\downarrow & & \downarrow & & \\ 1 & \to & \tilde{G} & \to & \tilde{G}.S & \to & S/F & \to & 1 \end{array}$$

and also the following commutative diagram

$$\begin{array}{ccccc} G(k') & \xrightarrow{\beta} & \bar{G}(k') \times (S/F)(k') & \xrightarrow{\delta'} & \mathrm{H}^1(k', \tilde{F}) \\ \downarrow & & \downarrow p' & & \downarrow q' \\ G(k') & \xrightarrow{\alpha} & (S/F)(k') & \xrightarrow{\delta'} & \mathrm{H}^1(k', \tilde{G}) \end{array}$$

By our assumption the Corestriction Principle holds for the image of $\delta'$. We claim that the composition of the maps

$$\bar{G}(k') \times (S/F)(k') \xrightarrow{p'} (S/F)(k') \xrightarrow{Cores_{S/F}} (S/F)(k) \xrightarrow{\delta} \mathrm{H}^1(k, \tilde{G}),$$

and that of the maps

$$\bar{G}(k') \times (S/F)(k') \xrightarrow{\delta'} \mathrm{H}^1(k', \tilde{F}) \xrightarrow{Cores_{\tilde{F}}} \mathrm{H}^1(k, \tilde{F}) \xrightarrow{q} \mathrm{H}^1(k, \tilde{G})$$

are the same. Indeed, denote by $p$ and $q'$ the maps similar to $p'$ and $q$, by considering the fields $k$ and $k'$ interchanged. Then for $x = (g', s') \in$

$\bar{G}(k') \times (S/F)(k')$ and $s = Cores_{S/F}(s') \in (S/F)(k)$ we have

$$\delta(Cores_{S/F}(p'(x))) = \delta(Cores_{S/F}(s')) = \delta(s).$$

By assumption there is $g \in \bar{G}(k)$ such that

$$Cores_{\tilde{F}}(\delta'(g')) = \delta(g),$$

hence

$$Cores_{\tilde{F}}(\delta'(g', s')) = \delta(g, s).$$

Since $p$ and $p'$ are surjective and the above diagram is commutative, it follows that for $y = (g, s) \in \bar{G}(k) \times (S/F)(k)$ we have

$$\begin{aligned}
\delta(Cores_{S/F}(p'(x)) &= \delta(s) \\
&= \delta(p(g, s)) \\
&= q(\delta(g, s)) \\
&= q(Cores_{\tilde{F}}(\delta'(g', s'))) = q(Cores_{\tilde{F}}(\delta'(x)))
\end{aligned}$$

as claimed. Now the assertion of the theorem follows from the equality $\alpha = p'\beta$. Indeed, let $x' \in G(k')$, $x'' = \beta(x') = (g', s')$, $y' = \alpha(x')$, $y = Cores_{S/F}(y')$. Then $\alpha(x') = p'\beta(x') = p'(x'')$ hence

$$\begin{aligned}
\delta(Cores_{S/F}(p'(x''))) &= q(Cores_{\tilde{F}}(\delta'(x''))) \\
&= q(Cores_{\tilde{F}}(\delta'(\beta(x')))) \\
&= 1,
\end{aligned}$$

since $\delta'\beta = 0$. Therefore

$$Cores_{S/F}(p'(x'')) \in \text{Ker } (\delta) = \text{Im } (\alpha).$$

Now consider the case $p = 1$. We consider the diagrams of cohomologies derived from the above diagram over $k$ and over $k'$. By Lemma 2.3 for any extension $K$ of $k$ there exist a pointed set, denoted by $\mathbf{H^2}(K)$ with morphisms of pointed sets such that the following diagram is commutative with

exact lines :

$$\begin{array}{ccccc} \mathrm{H}^1(k',G) & \overset{\beta_*}{\to} & \mathrm{H}^1(k',\bar{G}\times(S/F)) & \to & \mathrm{H}^2(k',\tilde{F}) \\ \downarrow & & \downarrow & & \downarrow \\ \mathrm{H}^1(k',G) & \overset{\alpha}{\to} & \mathrm{H}^1(k',S/F) & \overset{\delta''}{\to} & \mathbf{H^2}(k'). \end{array}$$

The meaning of introducing the set $\mathbf{H^2}(.)$ is to replace some $\mathrm{H}^2$-cohomology sets, which behave non-functorially (see e.g. [Gi], [Sp]), by some "cohomology" pointed set which makes our diagrams commutative. It is possible indeed, because sometimes we just treat cohomology sets as a "local" objects which make our diagram commutative as desired. So thinking of $\mathbf{H^2}$ as a "cohomology of something" (which exists as we have proved before) we may immitate the arguments for the case $p=0$ above. ∎

The following is in a sense a converse statement of what we have proved above.

**2.8. Proposition.** *Assume that the Corestriction Principle for images holds for* $\mathrm{H}^p(k,G)\to\mathrm{H}^p(k,T)$ *for all G, T as above with p=0 (resp. p=1). Then the same holds for* $\mathrm{H}^p(k,\bar{G})\to\mathrm{H}^{p+1}(k,F')$ *for any G, F' as above with p=0 (resp. p=1). In particular, if* a) *holds then* b) *holds.*

*Proof.* We consider in fact a slightly more general situation. Let us be given any isogeny $1\to F\to G_1\to G\to 1$ of connected reductive $k$-groups, with $F$ finite central $k$-subgroup of $G_1$ of multiplicative type. We will prove that the assumption of the proposition implies that the Corestriction Principle for images holds for $\mathrm{H}^p(k,G)\to\mathrm{H}^{p+1}(k,F)$, $p=0,1$.

To prove the assertion, we use the Ono's crossed diagram (see [O] for details) which allows one to embed an exact sequence with finite kernel of multiplicative type (i.e. isogeny) into another with quasi-split torus as a kernel. We will denote all maps in the following diagrams (for the level $k$ and $k'$) by the same symbols :

$$
\begin{array}{ccccccccc}
 & & 1 & & 1 & & & & \\
 & & \downarrow & & \downarrow & & & & \\
1 & \rightarrow & F & \rightarrow & G_1 & \stackrel{\alpha}{\rightarrow} & G & \rightarrow & 1 \\
 & & \downarrow & & \downarrow & & \downarrow = & & \\
1 & \rightarrow & T_1 & \rightarrow & H & \stackrel{\alpha}{\rightarrow} & G & \rightarrow & 1 \\
 & & \gamma\downarrow & & \downarrow\gamma & & & & \\
 & & T & = & T & & & & \\
 & & \downarrow & & \downarrow & & & & \\
 & & 1 & & 1 & & & &
\end{array}
$$

where $T_1$ is a quasi-split torus. From this diagram we derive the following commutative diagram

$$
\begin{array}{ccccccc}
 & 1 & & 1 & & & T(k') \\
 & \downarrow & & \downarrow & & & \downarrow \delta \\
1 \to & F(k') & \to & G_1(k') & \xrightarrow{\lambda} & G(k') & \xrightarrow{\beta} \mathrm{H}^1(k',F) \\
 & \downarrow & & \downarrow & & \downarrow = & \downarrow \theta \\
1 \to & T_1(k') & \to & H(k') & \xrightarrow{\alpha} & G(k') & \to 1 \\
 & \downarrow \pi & & \downarrow \gamma & & & \\
 & T(k') & = & T(k') & & & \\
 & \downarrow \delta & & \downarrow \zeta & & & \\
G(k') \xrightarrow{\beta} & \mathrm{H}^1(k',F) & \xrightarrow{i_*} & \mathrm{H}^1(k',G_1). & & &
\end{array}
$$

We need the following simple lemmas.

**2.8.1. Lemma.** [M1] *We have the following anti-commutative diagram*

$$
\begin{array}{ccc}
H(k') & \xrightarrow{\alpha} & G(k') \\
\gamma \downarrow & & \downarrow \\
T(k') & \xrightarrow{\delta} & \mathrm{H}^1(k',F)
\end{array}
$$

*for all field extension $k \subset k'$.*

We continue the proof of 2.8 and we assume first that $p = 0$.

Since $T_1$ is quasi-split, we have $\alpha(H(k')) = G(k')$. Now for any $g' \in G(k')$, let $h' \in H(k')$ such that $\alpha(h') = g'$, and denote $t' = \gamma(h')$, $f' = \beta(g')$. Since the diagram in Lemma 2.8.1 is anti-commutative, we have

$$
\begin{aligned}
\delta(t') &= \delta(\gamma(h')) \\
&= -\beta(\alpha(h'))
\end{aligned}
$$

$$= -f'.$$

Then $\theta(y') = 0$ hence $y' = \delta(z')$, $z' \in T(k')$. The image $z \in T(k)$ of $z'$ via $N_{k'/k} : T(k') \to T(k)$ is such that $\delta(z) = y := Cores_F(y')$. Now look at the diagram on the left hand side. Since the Corestriction Principle holds for the image of $H(k) \to T(k)$, there is $h \in H(k)$ such that

$$\gamma(h) = t := N_{k'/k}(t').$$

Let $g = \alpha(h)$. Then

$$\begin{aligned} \delta(t) &= \delta(\gamma(h)) \\ &= -\beta(\gamma(h) \\ &= -\beta(g) \\ &= \delta(N_{k/k}(t')) \\ &= Cores_F(\delta(t')) \\ &= -f, \end{aligned}$$

so $\beta(g) = f$ and $f \in \mathrm{Im}\ (\beta)$ and the case $p = 0$ is proved.

Now let $p = 1$. For any finite extension $k'$ of $k$, $g'$ any element from $\mathrm{H}^1(k', G)$ we choose a $g'$-lifting $z$-extension

$$1 \to T_1 \to H \to G \to 1,$$

defined over $k$, such that there is an embeding $F \hookrightarrow T_1$. We consider the following diagram, which is similar to the one we have just considered, with the only difference that the dimension is shifted.

$$
\begin{array}{ccccccc}
 & & & & & & \mathrm{H}^1(k',T) \\
 & & & & & & \downarrow \delta \\
\mathrm{H}^1(k',F) \to & \mathrm{H}^1(k',G_1) & \xrightarrow{\lambda} & \mathrm{H}^1(k',G) & \xrightarrow{\beta} & & \mathrm{H}^2(k',F) \\
\downarrow & \downarrow & & \downarrow = & & & \downarrow \theta \\
\mathrm{H}^1(k',T_1) \to & \mathrm{H}^1(k',H) & \xrightarrow{\alpha} & \mathrm{H}^1(k',G) & \to & & \mathrm{H}^2(k',T_1) \\
\downarrow \pi & \downarrow \gamma & & & & & \\
\mathrm{H}^1(k',T) = & \mathrm{H}^1(k',T) & & & & & \\
\downarrow \delta & & & & & & \\
\end{array}
$$

$$\mathrm{H}^1(k',G) \xrightarrow{\beta} \mathrm{H}^2(k',F)$$

We need the following analog of 2.8.1 for higher dimension

**2.8.2. Lemma.** *We have the following anti-commutative diagram*

$$
\begin{array}{ccc}
\mathrm{H}^1(k',H) & \xrightarrow{\alpha} & \mathrm{H}^1(k',G) \\
\Big\downarrow\gamma & & \Big\downarrow\beta \\
\mathrm{H}^1(k',T) & \xrightarrow{\Delta} & \mathrm{H}^2(k',F)
\end{array}
$$

*Proof.* Let $h = [(h_s)] \in \mathrm{H}^1(k',H),\ g = \alpha(h) \in \mathrm{H}^1(k',G)$. Then $g = [(g_s)]$, where $g_s = \alpha(h_s)$. We choose for each $s$ an element $g'_s \in G_1(k_s)$ such that $g_s = \alpha(g'_s)$. Then

$$h_s = g'_s t_s,\ t_s \in T_1(k_s).$$

One deduces from this

$$f_{sr}g'_{sr} = g'_s \, {}^s g'_r \ (s, r \in Gal(k_s/k'))$$

for some $f_{sr} \in F(k_s)$ and we know (see [Se]) that $(f_{sr})$ is 2-cocycle which is a representative of $\beta([(g_s)])$. From $h_s = g'_s t_s$ we deduce

$$\gamma(h_s) = \gamma(g'_s)\gamma(t_s) = \gamma(t_s),$$

hence for $t = [\gamma(h_s)] \in \mathrm{H}^1(k', T)$ we have

$$\Delta(t) = [(t_{sr}^{-1} t_s \, {}^s t_r)] \in \mathrm{H}^2(k', F).$$

Now the product of two 2-cocycles is

$$(t_{sr}^{-1} t_s \, {}^s t_r)(g'^{-1}_{sr} g'_s \, {}^s g'_r) = h_{sr}^{-1} h_s \, {}^s h_r = 1,$$

since $(h_s)$ is a 1-cocycle. Thus

$(*) \beta(\alpha(h)) = -\Delta(\gamma(h)),$

and the lemma follows. ∎

With $g' \in \mathrm{H}^1(k', G)$ as above let $h' \in \mathrm{H}^1(k', H)$ such that $g' = \alpha(h')$. Take a cocycle representative $(g_s)$ of $g'$ and let $g_s = \alpha(g_{1,s}), g_{1,s} \in G_1(k_s)$. Let $(h'_s)_s$ be a representative of $h'$, $h'_s \in H(k_s)$. Then

$$\beta(g') = [(g_{1,st}^{-1} g_{1,s} \, {}^s g_{1,t})] = -\Delta(\gamma(h'))$$

by the lemma above. Therefore

$$Cores_F(\beta(g')) = -Cores_F(\Delta(\gamma(h'))$$

$$= -\Delta(Cores_T(\gamma(h')).$$

By assumption, we have $Cores_T(\gamma(h')) = \gamma(h)$ for some $h \in \mathrm{H}^1(k, H)$. Let $g = \alpha(h) \in \mathrm{H}^1(k, G)$. Then

$$Cores_F(\beta(g')) = -\Delta(\gamma(h))$$

$$= \beta(\alpha(h)) \text{ (by } (*))$$

$$= \beta(g)$$

as required. ∎

Finally by summing up the results we proved above we obtain the following theorem which is the main result of this section. For the statements $a)$ - $d)$ considered above, let us denote by $x(p,q)$ the statement $x)$ evaluated at $(p,q)$, for $0 \leq p \leq q \leq 2$. For example, $a(1,2)$ means the statement $a)$ with $p=1, q=2$.

**2.9. Theorem.** 1) *All statements* $a)$ - $d)$ *are equivalent.*

2) *We have the following interdependence between the statements* $a)$ - $d)$ *with particular values of* $p$ *and* $q$.
$a)$ *For lower dimension :*

$$\begin{array}{c} a(0,0) \Leftrightarrow b(0) \Leftrightarrow c(0) \Leftrightarrow d(0) \\ \Downarrow \\ a(0,1) \end{array}$$

$b)$ *For higher dimension :*

$$\begin{array}{c} a(1,1) \Leftrightarrow b(1) \Leftrightarrow c(1) \Leftrightarrow d(1) \\ \Downarrow \\ a(1,2) \end{array}$$

*where two statements in the same row are connected by* $\Leftrightarrow$ *if they are equivalent and the down arrow indicates that the statements standing above imply the ones standing below.*

*Proof.* We just indicate the logical dependence of the statements of 1).
$d) \Rightarrow a)$ : see 2.1.
$a) \Rightarrow d)$ : see 2.4.
$b) \Leftrightarrow c)$ : see 2.5.
$c) \Rightarrow a)$ : see 2.7.
$a) \Rightarrow b)$ : see 2.8. ∎

From the proofs of propositions above we derive several consequences.

**2.10. Corollary.** *If either one of the conditions* $a)$ *or* $d)$ *holds (e.g. if*

*k ia a local or global field of characteristic 0) then for any isogeny of connected reductive k-groups*

$$1 \to F \to G_1 \to G_2 \to 1,$$

*the Corestriction Principle for images holds for standards maps*

$$\mathrm{H}^p(k, G_2) \to \mathrm{H}^{p+1}(k, F),\ p = 0, 1.$$

**2.11. Remarks.** 1) From the proof of Theorem 1.2, its corollary and Theorem 2.9 one may deduce the following new proof of Deligne's result mentioned above ([De, Prop. 2.4.8]) in the case $k$ is a local or global field of characteristic 0. Another proof is due to Milne and Shih [MS, Section 3]. (Unfortunately this proof is not a short one.)

**Corollary.** *If $k$ is a local or global field of characteristic 0 then $d(0)$ holds. In particular Theorem 1.1 holds.*

*Proof.* By Theorem 2.9, $d(0) \Leftrightarrow c(0)$. The proof of Theorem 1.2 reduces the proof of $c(0)$ to proving the Corestriction Principle for the kernel of $\mathrm{H}^1(k, \tilde{F}) \to \mathrm{H}^1(k, \tilde{G})$ so to the same thing for $\mathrm{H}^1(k, T) \to \mathrm{H}^1(k, \tilde{G})$ for a maximal $k$-torus $T$ which was what Theorem 1.3 asserted. ∎

2) A known sufficient condition for $c(0)$ to hold is that the group of $R$-equivalence of $G$ over $k'$ is trivial, i.e., $G(k')/R = 1$, since the Norm Theorem for the group of elements $R$-equivalent to 1 holds (see [G1, Prop. 3.3.2]). In [M1, Theorem 1] Merkurjev proved, among other results, a Norm Theorem from which the above result of [G1] follows. In Section 2 below we give a new proof of this result of Merkurjev.

3) The proof of Theorem 2.9 reduces the proof of Corestriction Principle for images for connected reductive groups to that of the maps

$$\mathrm{H}^p(k, \bar{G}) \to \mathrm{H}^{p+1}(k, \tilde{F}),$$

where $\tilde{F}$ is the center of a simply connected semisimple $k$-group $\tilde{G}$ with adjoint group $\bar{G}$. It is clear that we can reduce further to the case where $\bar{G}$ is almost simple. In this case, the Corestriction Principle for images is known for the case ${}^1\mathrm{A}_n$, $\mathrm{B}_n$ (due to the rationality of $\bar{G}$ and the result of Gille and Merkurjev mentioned above), $\mathrm{C}_n$ (due to Example 2 in Introduction above). It seems possible that it is also true for the case $\mathrm{D}_n$, since we see from above that the Corestriction Principle for images holds for the maps $G(k) \to \mathrm{H}^1(k, \mu_2)$, where $\Phi$ is a form of type $\mathrm{D}_n$, $G = \mathrm{SU}(\Phi)$ or its adjoint group. In general, according to Merkurjev [M2], the adjoint groups with non-trivial $R$-equivalence groups, hence *non* stably rational (even over number fields), exist.

# 3 Corestriction Principle for R-equivalence groups.

Let $G$ be a $k$-group. Two points $x$, $y \in G(k)$ are called $R-equivalent$ (after Manin) if there is a map $f : \mathbf{P}^1 \to G$ defined over $k$ and regular at 0 and 1, such that $f(0) = x$ and $f(1) = y$ (see [CTS] for more details). The subset $R := RG(k)$ of all elements of $G(k)$ which are $R$-equivalent to the identity is a *normal subgroup* of $G(k)$. It is well-known (see [CTS]) that for a field $k$ of characteristic 0, the factor group $G(k)/R$, called *the group of R-equivalence classes* of $G$ over $k$, is a birational invariant of the group $G$. In general, the study of the group $G(k)/R$ provides interesting information about the arithmetico-group-theoretic structure of the group $G(k)$, especially because there are many (even semisimple) groups with non-trivial R-equivalence groups (even over number fields).

In this section we are interested in the Corestriction Principle for images for $G(k)/R$ over local and global fields of characteristic 0. In [G2] it has been shown that for any reductive (hence also any) group $G$ defined over a number field $k$, the group of $R$-equivalences of $G$ over $k$ is finite. We use the notion of standard maps introduced in I. In [G1] Gille proved the following

**3.1. Theorem.** [G1] *Let $\pi : \tilde{G} \to G$ be an isogeny of connected reductive groups, all defined over a field $k$ of characteristic 0. Let $F = \mathrm{Ker}\,\pi$. Then for any finite extension $k'/k$ and the coboundary maps $\delta' : G(k') \to \mathrm{H}^1(k', F), \delta : G(k) \to \mathrm{H}^1(k, F)$ we have*

$$Cores_{k'/k}((\delta'(G(k'))) \subset \delta(G(k)).$$

In [M1] Merkurjev deduced 3.1 from the following result.

**3.2. Theorem.** [M1] *Let $\pi : G \to T$ be a homomorphism of connected reductive groups, where $T$ is a torus, all defined over a field $k$ of characteristic 0. Then*

$$N_{k'/k}(\pi(RG(k')) \subset \pi(RG(k)).$$

We give here a new proof of 3.2 by using 3.1 and the reductions made in Section 2.

*Proof.* Let $G' = [G, G]$, $T' = G/G'$. It is obvious that if 3.2 is true for the pair $(G, T')$ then it is true for $(G, T)$. So we may assume that $T = G/G'$. By Lemma 1.6 there exists a $z$-extension $\pi_1 : G_1 \to T_1$ of $\pi : G \to T$. If we denote by $\alpha : G_1 \to G$, $\beta : T_1 \to T$ the corresponding projections, then it is easy to see that (see e.g. [T2])

$$\alpha(RG_1(K)) = RG(K), \beta(RT_1(K)) = RT(K)$$

for any extension $K/k$. Therefore we may assume that $G$ has simply connected semisimple part. It is obvious that if the lemma is true for some power $G^n = G \times \cdots \times G$, then it is also true for $G$, so by virtue of Lemma 1.10 of [S] (used before, in the proof of Theorem 1.3) we may assume that $G$ has a special covering $\tilde{G} \times T'$, where $T'$ is an induced $k$-torus, and we have the following exacts sequence of algebraic groups, all defined over $k$ :

$$1 \to F \to \tilde{G} \times T' \to G \to 1,$$

where $F$ is a finite central subgroup of $\tilde{G} \times T'$. From this we derive the following $3 \times 3$-commutative diagram

$$\begin{array}{ccccccccc}
 & & 1 & & 1 & & 1 & & \\
 & & \downarrow & & \downarrow & & \downarrow & & \\
1 & \to & F \cap \tilde{G} & \to & F & \stackrel{u}{\to} & F' & \to & 1 \\
 & & \downarrow & & \downarrow & & \downarrow & & \\
1 & \to & \tilde{G} & \to & \tilde{G} \times T' & \to & T' & \to & 1 \\
 & & l\downarrow & & \downarrow\pi & & \downarrow & & \\
1 & \to & \tilde{G} & \to & G & \to & T & \to & 1 \\
 & & \downarrow & & \downarrow & & \downarrow & & \\
 & & 1 & & 1 & & 1 & &
\end{array}$$

Since $\tilde{G}$ is simply connected, $l$ is an isomorphism, hence $F \cap \tilde{G} = 1$, and $u$ is also an isomorphism. Denote by $\mathrm{RH}^1(K, F)$ the set of all elements R-equivalent to the trivial element of $\mathrm{H}^1(K, F)$. From the diagram above we derive the following commutative diagram

$$\begin{array}{ccccc}
R(\tilde{G}(k) \times T'(k)) & \stackrel{\pi}{\to} & RG(k) & \stackrel{\delta_G}{\to} & \mathrm{RH}^1(k, F) \\
\downarrow p & & \downarrow q & & \simeq \downarrow r \\
RT'(k) & \to & RT(k) & \stackrel{\delta_T}{\to} & \mathrm{RH}^1(k, F')
\end{array}$$

where $r$ is an isomorphism, induced from $u$. We have similar diagram when $k'$ replaces $k$, where one changes $p \to p'$, etc..., for example

$$\delta'_G : G(k') \to \mathrm{H}^1(k', F), \delta'_T : T(k') \to \mathrm{H}^1(k', F)$$

are coboundary maps. Let $g' \in G(k')$. If $g' \in RG(k')$, then we have

$$\delta_T(N_{k'/k}(q'(g'))) = N_{k'/k}(\delta'_T(q'(g')))$$

$$= N_{k'/k}(r'(\delta'_G(g')))$$

$$= r(N_{k'/k}(\delta'_G(g'))).$$

By assumption, $N_{k'/k}(\delta'_G(g')) = \delta_G(h)$ for some $g \in RG(k)$. Therefore

$$\delta_T(N_{k'/k}(q'(g'))) = r(\delta_G(g)) = \delta_T(q(g)),$$

thus

$$N_{k'/k}(q'(g')) = q(g)t,$$

where $t \in \text{Ker }(\delta_T) = \text{Im }(T'(k) \to T(k))$. Since $T'$ is an induced $k$-torus, $t \in RT(k)$, and since $p$ is just the projection, from the above commutative diagram we deduce $t \in q(RG(k))$. Thus

$$N_{k'/k}(q'(g')) \in q(RG(k))$$

as required. ∎

We derive the following consequence.

**3.3. Theorem.** *Assume that $k$ is a local or global field of characteristic 0. Then for any connected reductif $k$-group $G$, a $k$-torus $T$, a standard map $\pi : G(k) \to T(k)$ and for any finite extension $k'$ of $k$, the norm homomorphism $T(k') \to T(k)$ induces a canonical functorial norm map for images*

$$N_{k'/k} : \text{Im }(G(k')/R \to T(k')/R) \to \text{Im }(G(k)/R \to T(k)/R).$$

*Proof.* . We know from Section 1 that

$$N_{k'/k}(\text{Im }\ (G(k') \to T(k')) \subset (\text{Im }\ (G(k) \to T(k)).$$

By Theorem 3.2

$$N_{k'/k}(\text{Im }(RG(k') \to RT(k')) \subset (\text{Im }(RG(k) \to RT(k)),$$

and the theorem follows from these two inclusions.

Another proof is as follows. First one reduces (as in Section 2) the proof to the case where $G$ has simply connected semisimple part $\mathcal{D}G = [G, G]$ and $T = G/\mathcal{D}G$ is the torus quotient of $G$. Then the theorem follows from the fact that the natural projection $G \to T$ induces a *surjective* map $G(k)/R \to T(k)/R$ by Theorem 4.12 of [T3]. ∎

**3.4. Corollary**. *With above notation, for any isogeny of connected k-groups*

$$1 \to F \to H \to G \to 1,$$

*with finite F, the Corestriction Principle for images holds for the map*

$$G(k)/R \to (\mathrm{Im}\ (\delta))/R,$$

*where* $\delta$ *is the connecting map* $G(k) \to \mathrm{H}^1(k, F)$*, and the R-equivalence in* $Im(\delta)$ *is induced from that of* $G(k)$ *as defined in [G1].*

*Proof.* Use the same Ono's crossed diagram as in the proof of Proposition 2.8 in Section 2. ∎

# 4 Knebusch Norm Principle

Let $k$ be a field of characteristic $\neq 2$, $q$ a non-degenerate quadratic form over $k$, and $K$ any finite extension of degree $n$ of $k$. The Knebusch Norm Principle (see [L]) states that for any $x \in D(q \otimes K)$, the set of values of the form $q$ in $K^*$, $N_{K/k}(x)$ is the product of $n$ elements from $D(q)$. In particular, for the group $D[q] := \langle D(q) \rangle$ generated by the non-zero values of $q$ the Norm Principle holds :

$$N_{K/k}(D[q \otimes K]) \subset D[q].$$

The first natural question arises :

**Question 1.** *What happens if the quadratic form is replaced by a homogeneous form of degree* $\geq 3$ *?*

In this section we are interested in the following other natural questions. Let $\Phi$ be a hermitian form as in Example 2 of Introduction. Denote by $D(\Phi)$ (resp. $D[\Phi]$) the set of (resp. the group generated by) non-zero values of $\Phi$. Let $ND(\Phi)$ (resp. $ND[\Phi]$) the image of $D(\Phi)$ (resp. $D[\Phi]$) in $k^*$ via $Nrd_{D/k_0}$. For any finite extension $K/k$ of degree $n$ we ask

**Question 2.** *Is there any function $f(n)$ with values in $\mathbf{N}$ such that for any $x \in ND(\Phi \otimes K)$ , $N_{K/k}(x)$ is the product of $f(n)$ elements from $ND(\Phi)$ ?*

**Question 3.** *When does $N_{K/k}(ND[\Phi \otimes K]) \subset ND[\Phi]$ ?*

Concerning the hermitian forms, we assume $k$ is a local or global field of characteristic $\neq 2$. To answer the questions above, it is easy to see that, due to Knebusch's Norm Theorem above, it is sufficient to consider the cases where $k$ is a global field and $\Phi$ is of type A or D and the division $k$-algebra $D$ is non-trivial. We will see that the validity of the Norm Principle depends very much on the arithmetic nature of the base field.

From now on we assume that $k$ is a global field of characteristic $\neq 2$ and that $D$ is non trivial. First assume that $\Phi$ is of type A. For a valuation $v$ of $k$ we denote by $k_v$ the completion of $k$ at $v$. Let $n = \text{rank}(\Phi)$.

Let $n = 1$, $\Phi = \langle d \rangle$, $d^J = d$. Then

$$D(\Phi) = \{x^J dx : x \in D^*\}.$$

Hence

$$ND(\Phi) = \{Nrd_{D/k_0}(d) Nrd_{D/k_0}(x)^J Nrd_{D/k_0}(x) : x \in D^*\}$$

and for finite extension $K/k$, $K_0 = Kk_0$ and $z \in ND(\Phi \otimes K)$ we have $z = Nrd_{D_K/K_0}(d) Nrd_{D_K/K_0}(x)^J Nrd_{D_K/K_0}(x)$ for $x \in D \otimes K$. Since

$$N_{K_0/k_0}(Nrd_{D_K/K_0}(D_K)) \subset Nrd_{D/K}(D),$$

the Norm Principle is true.

Now let $n \geq 2$. Since the Strong Hasse Principle holds for $\Phi$ (see [Sc, Chapter 10]), and $\Phi_v := \Phi \otimes k_v$ is equivalent to a hermitian form over $k_v$ for all nonarchimedean local fields $k_v$, we have

$$\begin{aligned} D(\Phi) &= \{x \in D^* : x^J = x, \Phi_x := \langle -x \rangle \perp \Phi \text{ is isotropic } \} \\ &= \{x \in D^* : x^J = x, \Phi_{x,v} \text{ is isotropic for all real places } v \in V_0\} \end{aligned}$$

where $V_0 = \{v \in \infty : \Phi_v$ is anisotropic $\}$.

Denote by $d$ the degree of $D$, $s_v$ the signature of $\Phi_v$ for $v \in V_0$. By definition of $V_0$, $s_v = nd$ (resp. $-nd$) if $\Phi_v$ is positively (resp. negatively) definite. Let $V_0^+ = \{v \in V_0 : s_v > 0\}$ and $V_0^- := V_0 \setminus V_0^+$. Hence

$$D(\Phi) = \{x \in D^* : x^J = x,\ sign_v(\langle x \rangle_v) = \epsilon d, \forall v \in V_0^\epsilon, \epsilon = \pm\}.$$

Now from [Sc, Chapter 10, 6.8, 6.9 and 6.10] it follows that

$$ND(\Phi) = k(+,-),$$

where

$$k(+,-) = \{\delta \in k^* : \delta > 0\ (\forall v \in V_0^+),\ \eta\delta > 0\ (\forall v \in V_0^-)\},$$

where $\eta = (-1)^d$. In other words, $k(+,-)$ is the subset of all elements of $k^*$ which have certain assigned signs at $V_0$. Denote by $V_{0,K}$, $V_{0,K}^+$, $V_{0,K}^-$, $K(+,-)$ the similar sets when $k$ is replaced by $K$. We can check without difficulty whether $N_{K/k}(K(+,-)) \subset k(+,-)$. The case of forms of type D can be considered in a similar way, by making use of the Kneser's Strong Hasse Principle for skew-hermitian forms of dimension $\geq 3$.

**Acknowledgements.** I would like to thank M. Borovoi for sending preprints of his works, J. -L. Colliot-Thélène (while at the Fields Institute), M. Kolster, C. Riehm and V. Snaith for discussions, Department of Mathematics and Statistics, McMaster University, J. Sonn and Lady Davis Foundation for support during the preparation of this paper.

Some of the main results of this paper has been reported in the author's talk at Mathematisches Institut Oberwolfach (Germany) in June 1995. I wish to thank M. Knebusch, W. Scharlau and A. Pfister for the help and hospitality during my stay in Oberwolfach.